\documentclass[10pt]{article}
\usepackage[latin1]{inputenc}
\usepackage{epsfig}
\usepackage{color}
\usepackage[british,english]{babel}
\usepackage{amsthm}
\usepackage{amsmath}
\usepackage{amsfonts}
\usepackage{amssymb}
\usepackage{graphicx}

\usepackage{fancyhdr}

\newtheorem{corollary}{Corollary}[section]
\newtheorem{definition}[corollary]{Definition}

\newtheorem{lemma}[corollary]{Lemma}
\newtheorem{proposition}[corollary]{Proposition}

\newtheorem{theorem}[corollary]{Theorem}

\date{}
\begin{document}
\title{Reaction-diffusion equation on thin porous media}\maketitle

\vskip-30pt
 \centerline{Mar\'ia ANGUIANO\footnote{Departamento de An\'alisis Matem\'atico. Facultad de Matem\'aticas.
Universidad de Sevilla, 41012 Sevilla (Spain)
anguiano@us.es}}

 \renewcommand{\abstractname} {\bf Abstract}
\begin{abstract} 
We consider a reaction-diffusion equation on a 3D thin porous media of thickness $\varepsilon$ which is perforated by periodically distributed  cylinders of size $\varepsilon$. On the boundary of the cylinders we prescribe a dynamical boundary condition of {\it pure-reactive} type. As $\varepsilon\to 0$, in the 2D limit the resulting reaction-diffusion equation has a source term coming from the dynamical-type boundary conditions imposed on boundaries of the original 3D domain.
\end{abstract}

 {\small \bf AMS classification numbers:} 35K57; 35B27  
 
 
 \section {Introduction and setting of the problem}\label{S1} 
 We consider a parabolic problem in a thin porous media $\Omega_\varepsilon$ of thickness $\varepsilon$ which is perforated by periodically distributed cylinders (obstacles) of size $\varepsilon$. On the boundary of the cylinders, we prescribe a dynamical boundary condition of {\it pure-reactive} type. The aim of this work is to prove the convergence of the homogenization process when $\varepsilon$ goes to zero.
 
The asymptotic behavior of the solution of linear parabolic problems with dynamical boundary conditions of {\it pure-reactive} type in periodically perforated domains was analyzed by Timofte in \cite{Timofte}. The author obtained a new parabolic limiting problem defined on a unified domain with extra terms coming from the influence of the dynamical boundary conditions. More recently, in \cite{Anguiano_MJOM, Anguiano_ZAMM}, we study the nonlinear case where the nonlinearity appears reflected in the limit equation. However, to our knowledge, there does not seem to be in the literature any study of the asymptotic behavior of the solution of parabolic models associated with dynamical boundary conditions of {\it pure-reactive} type in a thin porous media. Let us introduce the model we will be involved with in this paper.

{\bf The geometrical setting.} The periodic porous media is defined by a domain $\omega$ and an associated microstructure, or periodic cell $Y^{\prime}=[-1/2,1/2]^2$, which is made of two complementary parts: the obstacle part $F^{\prime}$, and $Y'_{f}=Y'\setminus \overline F'$. More precisely, we assume that $\omega$ is a smooth, bounded, connected set in $\mathbb{R}^2$, and that $F^{\prime}$ is an open connected subset of $Y^\prime$ with a smooth boundary $\partial F^\prime$, such that $\overline F^\prime$ is strictly included  in $Y^\prime$.\\

The microscale of a porous media is a small positive number ${\varepsilon}$. The domain $\omega$ is covered by a regular mesh of square of size ${\varepsilon}$: for $k^{\prime}\in \mathbb{Z}^2$, each cell $Y^{\prime}_{k^{\prime},{\varepsilon}}={\varepsilon}k^{\prime}+{\varepsilon}Y^{\prime}$ is divided in an obstacle part $F^{\prime}_{k^{\prime},{\varepsilon}}$ and $Y^{\prime}_{f_{k^{\prime}},{\varepsilon}}$, i.e. is similar to the unit cell $Y^{\prime}$ rescaled to size ${\varepsilon}$. We define $Y=Y^{\prime}\times (0,1)\subset \mathbb{R}^3$, which is divided in an obstacle part $F=F'\times(0,1)$ and $Y_{f}=Y'_f\times (0,1)$, and consequently $Y_{k^{\prime},{\varepsilon}}=Y^{\prime}_{k^{\prime},{\varepsilon}}\times (0,1)\subset \mathbb{R}^3$, which is also divided in an obstacle part $F_{{k^{\prime}},{\varepsilon}}$ and $Y_{f_{k^{\prime}},{\varepsilon}}$. 

We denote by $\tau(\overline F'_{k',\varepsilon})$ the set of all translated images of $\overline F'_{k',\varepsilon}$. The set $\tau(\overline F'_{k',\varepsilon})$ represents the obstacles in $\mathbb{R}^2$.

We define $\omega_{\varepsilon}=\omega\backslash\bigcup_{k^{\prime}\in \mathcal{K}_{\varepsilon}} \overline F^{\prime}_{{k^{\prime}},{\varepsilon}}\subset \mathbb{R}^2,$ where $\mathcal{K}_{\varepsilon}=\{k^{\prime}\in \mathbb{Z}^2: Y^{\prime}_{k^{\prime}, {\varepsilon}} \cap \omega \neq \emptyset \}$, and we define the thin porous media $\Omega_{\varepsilon}\subset \mathbb{R}^3$ by 
\begin{equation}\label{Dominio1}
\Omega_{\varepsilon}=\{  (x_1,x_2,x_3)\in \omega_{\varepsilon}\times \mathbb{R}: 0<x_3<\varepsilon \}.
\end{equation}

We make the following assumption:
$$\text{The obstacles } \tau(\overline F'_{k',\varepsilon}) \text{ do not intersect the boundary } \partial \omega.$$

We define $F^\varepsilon_{k',\varepsilon}=F'_{k',\varepsilon}\times (0,\varepsilon)$. Denote by $S_\varepsilon$ the set of the obstacles contained in $\Omega_\varepsilon$. Then, $S_\varepsilon$ is a finite union of obstacles, i.e. $$S_\varepsilon=\bigcup_{k^{\prime}\in \mathcal{K}_{\varepsilon}} \overline F^\varepsilon_{k',\varepsilon}.$$

We define 
\begin{equation}\label{OmegaTilde}
\widetilde{\Omega}_{\varepsilon}=\omega_{\varepsilon}\times (0,1), \quad \Omega=\omega\times (0,1), \quad \Lambda_\varepsilon=\omega\times (0,\varepsilon).
\end{equation}

We observe that $\widetilde{\Omega}_{\varepsilon}=\Omega\backslash \bigcup_{k^{\prime}\in \mathcal{K}_{\varepsilon}} \overline F_{{k^{\prime}}, {\varepsilon}},$ and we define $F_\varepsilon=\bigcup_{k^{\prime}\in \mathcal{K}_{\varepsilon}} \overline F_{k^\prime, \varepsilon}$ as the set of obstacles contained in $\widetilde \Omega_\varepsilon$.\\
 
{\bf Position of the problem.} We consider the following problem for a reaction-diffusion equation with dynamical boundary conditions of {\it pure-reactive} type on the surface of the obstacles and zero Dirichlet condition on the exterior boundary,
\begin{equation}
\left\{
\begin{array}
[c]{r@{\;}c@{\;}ll}%
\displaystyle\frac{\partial u_\varepsilon}{\partial t}-\Delta\,u_\varepsilon+\kappa\, u_\varepsilon &
= &
f(x',t)\quad & \text{\ in }\;\Omega_\varepsilon\times(0,T) ,\\
\nabla u_\varepsilon \cdot \nu_\varepsilon+\varepsilon\,\displaystyle\frac{\partial u_\varepsilon}{\partial
t} & = & \varepsilon\,g(x',t) & \text{\ on }%
\;\partial S_\varepsilon\times( 0,T),\\
u_\varepsilon(x,0) & = & u_\varepsilon^{0}(x), & \text{\ for }\;x\in\Omega_\varepsilon,\\
u_\varepsilon(x,0) & = & \psi_\varepsilon^{0}(x), & \text{\ for
}\;x\in\partial S_\varepsilon,\\
u_\varepsilon&=& 0, & \text{\ on }%
\;\partial \Lambda_{\varepsilon}\times( 0,T),
\end{array}
\right. \label{PDE}%
\end{equation}
where $u_\varepsilon=u_\varepsilon(x,t)$, $x=(x',x_3)\in \Omega_\varepsilon$, $t\in (0,T)$ and $T>0$. The first equation states the law of standard diffusion in $\Omega_\varepsilon$, $\Delta=\Delta_x$ denotes the Laplacian operator with respect to the space variable and $\kappa>0$ is a given constant. 
The boundary equation (\ref{PDE})$_{2}$ is multiplied by $\varepsilon$ to compensate the growth of the surface by shrinking $\varepsilon$, where the value of $u_\varepsilon$ is assumed to be the trace of the function $u_\varepsilon$ defined for $x\in \Omega_\varepsilon$, $\nu_\varepsilon$ is the outer normal to $\partial S_\varepsilon$. The term $\nabla u_\varepsilon \cdot \nu_\varepsilon$ represents the heat flux across the boundary (see Goldstein \cite{Goldstein} for more details). We assume that
\begin{equation}\label{hyp 0}
u_\varepsilon^{0}\in L^2\left( \Omega\right),\quad
\psi_\varepsilon^{0}\in L^{2}\left( \partial S_\varepsilon\right),
\end{equation}
are given, and that there exists a positive constant $C$ such that
\begin{equation}\label{Initial_condition}
\varepsilon^{-1}|u_\varepsilon^0|^2_{\Omega_\varepsilon}+ |\psi_\varepsilon^0|^2_{\partial S_\varepsilon}\leq C,
\end{equation}
where $|\cdot|_{\Omega_\varepsilon}$ and $|\cdot|_{\partial S_\varepsilon}$ denote the norm in $L^2(\Omega_\varepsilon)$ and $L^2(\partial S_\varepsilon)$, respectively.

Moreover, we assume that
\begin{equation}\label{hyp 0'}
f\in L^{2}\left(0,T
;L^{2}\left( \omega\right) \right),\quad g\in
L^{2}\left(0,T;H_0^{1}\left(
\omega\right) \right),
\end{equation}
are given. This choice of $f$ and $g$ is usual when dealing with thin domains. Since the thickness  of the domain, $\varepsilon$,  is small then the exterior force can be considered independent of the vertical variable.

In this paper, our main motivation is to study the asymptotic behavior of the solution $u_{\varepsilon}$ of (\ref{PDE}) as $\varepsilon$, the size of the cylinders and the corresponding thickness of the thin film, tend to zero. We combine the homogenization techniques used in porous media and in thin film in order to rigorously derive the homogenized model.

For this purpose, we first perform a change of variables which consists in stretching in the $x_3$-direction by a factor $1/\varepsilon$. As in \cite{Anguiano_SuarezGrau}, we use the dilatation in the variable $x_3$, i.e.
\begin{equation}\label{dilatacion}
y_3=\frac{x_3}{\varepsilon},
\end{equation}
in order to have the functions defined in the open set with fixed height $\widetilde \Omega_\varepsilon$ defined by (\ref{OmegaTilde}). In this sense, we define $\tilde{u}_{\varepsilon}$ by 
\begin{eqnarray}\label{definition_u_tilde}
\tilde{u}_{\varepsilon}(x^{\prime},y_3,t)=u_{\varepsilon}(x^{\prime},\varepsilon y_3,t), \text{\ \ } a.e.\text{\ } (x^{\prime},y_3,t)\in \widetilde{\Omega}_{\varepsilon}\times (0,T),
\end{eqnarray}
and our goal then is to describe the asymptotic behavior of this new sequence $\tilde u_\varepsilon$ when $\varepsilon$ tends to zero.

Using the dilatation in the variable $x_3$ given by (\ref{dilatacion}), we have that $\Delta u_\varepsilon$ is transformed into
$$
\Delta_{x'}\tilde u_\varepsilon+{1\over \varepsilon^2}\partial^2_{y_3}\tilde u_\varepsilon.
$$
The factor $1/\varepsilon^2$ in front of the derivate in the $y_3$ direction means a very fast diffusion in the vertical direction. In some sense, we have substituted the thin porous media $\Omega_\varepsilon$ by a non-thin domain $\widetilde \Omega_\varepsilon$, but with a very strong diffusion mechanism in the $y_3$-direction. Because of the presence of this very strong diffusion mechanism, it is expected that solutions of the problem which satisfies $\tilde u_\varepsilon$, defined in (\ref{PDE_dilatado}), become homogeneous in the $y_3$-direction so that the limit solution will not have a dependence in this direction, and therefore, the limit problem will be two dimensional. 

The method we follow in this paper is the so-called energy method of Tartar \cite{Tartar0}, which has been considered by many authors (see, for instance, Cioranescu and Donato \cite{Ciora2}, Conca and Donato \cite{Conca1}, Conca {\it et al.} \cite{Conca2}, Timofte \cite{Timofte} and \cite{Anguiano_MJOM}) and the technique introduced by Vanninathan \cite{Vanni} for the Steklov problems which transforms surface integrals into volume integrals. We get the following limit problem to (\ref{PDE_dilatado}), as $\varepsilon$ goes to zero.

\begin{theorem}[Main Theorem]\label{Main}
Under the assumptions (\ref{Initial_condition})--(\ref{hyp 0'}), assume that $(u_\varepsilon^0,\psi_\varepsilon^0)\in H^1(\Omega_\varepsilon)\times H^{1/2}(\partial S_\varepsilon)$ and there exists a positive constant $C$ such that
\begin{equation}\label{Initial_condition_2}
\varepsilon^{-1/2}|\nabla u_\varepsilon^0|_{\Omega_\varepsilon}\leq C.
\end{equation}
Let $(u_\varepsilon, \psi_\varepsilon)$ be the unique solution of the problem (\ref{PDE}), where $\psi_\varepsilon$ is the trace of the function $u_\varepsilon$ defined for $x\in \Omega_\varepsilon$ and a.e. $t\in (0,T]$. Then, there exists an extension $\widetilde \varPi_\varepsilon\tilde u_\varepsilon$ of $\tilde u_\varepsilon$, where $\tilde u_\varepsilon$ is the dilation of $u_\varepsilon$ given by (\ref{definition_u_tilde}), into all $\Omega\times (0,T)$, such that, as $\varepsilon\to 0$,
$$\widetilde \varPi_\varepsilon\tilde u_\varepsilon(t) \to u(t) \quad \text{strongly in } L^2(\Omega),\quad \forall t\in[0,T],$$
where $u\in L^2(0,T;H_0^1(\omega))$ is the unique solution of the following problem
 \begin{equation}\label{limit_problem}
\left\{
\begin{array}{l}
\displaystyle \left({|Y'_f|\over |Y'|}+{|\partial F'| \over |Y'|} \right)\displaystyle\frac{\partial u}{\partial t}-{\rm div}_{x'}\left(Q\nabla_{x'}u \right)+ {|Y'_f|\over |Y'|}\kappa\, u
=\displaystyle {|Y'_f|\over |Y'|}f+{|\partial F'| \over |Y'|}g,   \text{\ in }\;\omega\times(0,T) ,\\[2ex]
u(x',0)  =  u_{0}(x'),  \text{\ for }\;x'\in\omega,\\[2ex]
u= 0,  \text{\ on }
\;\partial \omega\times( 0,T).
\end{array}
\right.
\end{equation}
The homogenized matrix $Q=((q_{i,j}))$, $1\leq i,j\leq 2$, which is symmetric and positive-definite, is given by 
\begin{equation}\label{matrix}
q_{i,j}={1\over |Y'|}\int_{Y'_f}\left(e_i+\nabla_{y'} w_i \right)\cdot \left(e_j+\nabla_{y'} w_j \right)dy',
\end{equation}
where $w_i\in \mathbb{H}_{{\rm per}}\setminus \mathbb{R}$, $ i=1,2$, is the unique solution of the cell problem 
\begin{equation}\label{system_eta}
\left\{
\begin{array}{l}
\displaystyle -{\rm div}_{y'}\left(e_i+\nabla_{y'} w_i \right)=0,   \text{\ in }Y'_f,\\[2ex]
(e_i+\nabla_{y'} w_i)\cdot \nu' =0,  \text{\ on }\partial F',\\[2ex]
w_i  \text{\ is }Y'-\text{periodic}.
\end{array}
\right.
\end{equation}
Here, $e_i$ is the $i$ element of the canonical basis in $\mathbb{R}^2$, $\nu'$ is the outer normal to $\partial F'$, and $\mathbb{H}_{{\rm per}}$ is the space of functions from $ H^1(Y'_f)$ which are $Y'$-periodic.
\end{theorem}
The structure of the paper as follows. In Section \ref{S2}, we introduce some notations which are used in the paper. In Section \ref{S3}, we give a weak formulation of the problem, and establish the existence and uniqueness of solution. In Section \ref{S4}, we establish the problem which satisfies $\tilde u_\varepsilon$ given by (\ref{definition_u_tilde}). Some {\it a priori} estimates for $\tilde u_\varepsilon$ are rigorously obtained in Section \ref{S5}. A compactness result is addressed in Section \ref{S6}. Finally, the main goal of proving the asymptotic behavior of the solution $\tilde u_\varepsilon$ is achieved in Section \ref{S7}.

\section{Some notations}\label{S2}
Along this paper, the points $x\in\mathbb{R}^3$ will be decomposed as $x=(x^{\prime},x_3)$ with $x^{\prime}\in \mathbb{R}^2$, $x_3\in \mathbb{R}$. We also use the notation $x^{\prime}$ to denote a generic vector of $\mathbb{R}^2$.
We denote by $\chi_{\widetilde \Omega_\varepsilon}$ the characteristic function of the domain $\widetilde \Omega_\varepsilon$.

We denote by $(\cdot,\cdot) _{\Omega_\varepsilon}$ (respectively, $(
\cdot,\cdot)_{\partial S_\varepsilon}$) the inner product in
$L^{2}(\Omega_\varepsilon)$ (respectively, in $L^{2}(\partial S_\varepsilon)$),
and by $\left\vert \cdot\right\vert _{\Omega_\varepsilon}$
(respectively, $\left\vert \cdot\right\vert
_{\partial S_\varepsilon}$) the associated norm. We also denote $(\cdot,\cdot) _{\Omega_\varepsilon}$ the inner product in $(L^2(\Omega_\varepsilon))^3$.

We denote by $(\cdot,\cdot) _{\widetilde \Omega_\varepsilon}$ (respectively, $(
\cdot,\cdot)_{\partial F_\varepsilon}$) the inner product in
$L^{2}(\widetilde \Omega_\varepsilon)$ (respectively, in $L^{2}(\partial F_\varepsilon)$),
and by $\left\vert \cdot\right\vert _{\widetilde \Omega_\varepsilon}$
(respectively, $\left\vert \cdot\right\vert
_{\partial F_\varepsilon}$) the associated norm. We also denote $(\cdot,\cdot) _{\widetilde \Omega_\varepsilon}$ the inner product in $(L^2(\widetilde \Omega_\varepsilon))^N$, with $N=2,3$.

We denote by $(\cdot,\cdot) _{\Omega}$ the inner product in
$L^{2}(\Omega)$,
and by $\left\vert \cdot\right\vert _{\Omega}$ the associated norm. We also denote $(\cdot,\cdot) _{\Omega}$ the inner product in $(L^2(\Omega))^3$. By $\left\Vert \cdot\right\Vert _{ \Omega}$ we denote the norm in
$H^{1}( \Omega)$.

We denote by $(\cdot,\cdot) _{\omega_\varepsilon}$ the inner product in $L^2(\omega_\varepsilon)$ and by $|\cdot|_{\omega_\varepsilon}$ the associated norm. We also denote $(\cdot,\cdot) _{\omega_\varepsilon}$ the inner product in $(L^2(\omega_\varepsilon))^2$. By $||\cdot ||_{\omega_\varepsilon}$ we denote the norm in $H^1(\omega_\varepsilon)$, which is associated to the inner
product $$((\cdot,\cdot))_{ \omega_\varepsilon}:=
(\nabla_{x'}\cdot,\nabla_{x'}\cdot)_{ \omega_\varepsilon}+\left(\cdot,\cdot\right)
_{ \omega_\varepsilon}.$$

By $\left\Vert \cdot\right\Vert _{\widetilde \Omega_\varepsilon}$ we denote the norm in
$H^{1}(\widetilde \Omega_\varepsilon)$, which is associated to the inner
product $$((\cdot,\cdot))_{\widetilde \Omega_\varepsilon}:=
(\nabla_{x'}\cdot,\nabla_{x'}\cdot)_{\widetilde \Omega_\varepsilon}+(\partial_{y_3}\cdot,\partial_{y_3}\cdot)_{\widetilde \Omega_\varepsilon}+\left(\cdot,\cdot\right)
_{\widetilde \Omega_\varepsilon}.$$

By $||\cdot||_{\widetilde \Omega_\varepsilon,T}$ we denote the norm in $L^2(0,T;H^1(\widetilde \Omega_\varepsilon))$. 

By $|\cdot|_{\widetilde \Omega_\varepsilon,T}$ (respectively $|\cdot|_{\partial F_\varepsilon,T}$), we denote the norms in $L^2(0,T;L^2(\widetilde \Omega_\varepsilon))$ and $L^2(0,T;(L^2(\widetilde \Omega_\varepsilon))^3)$ (respectively $L^2(0,T;L^2(\partial F_\varepsilon))$). By $|\cdot|_{ \Omega,T}$, we denote the norms in $L^2(0,T;L^2( \Omega))$ and $L^2(0,T;(L^2( \Omega))^3)$.

We denote by $\gamma_{0}$ the trace operator $u\mapsto
u|_{\partial\Omega_\varepsilon}$. The trace operator belongs to
$\mathcal{L}(H^1(\Omega_\varepsilon), H^{1/2}(\partial\Omega_\varepsilon))$, and we will
use $\|\gamma_0\|$ to denote the norm of $\gamma_0$ in this space. Analogously, we denote by $\gamma_{0}$ the trace operator $\tilde u\mapsto
\tilde u|_{\partial\widetilde \Omega_\varepsilon}$.

We will use $\|\cdot\|_{\partial\widetilde \Omega_\varepsilon}$ to denote the
norm in $H^{1/2}(\partial\widetilde \Omega_\varepsilon),$ which is given by
$\|\tilde \phi\|_{\partial\widetilde \Omega_\varepsilon}=\inf\{\|\tilde v\|_{\widetilde \Omega_\varepsilon}:\;
\gamma_0(\tilde v)=\tilde \phi\}$. We recall that with this norm,
$H^{1/2}(\partial\widetilde \Omega_\varepsilon)$ is a Hilbert space. 

We denote by $H^r_{\partial \Lambda_\varepsilon}(\Omega_\varepsilon)$ and $H^r_{\partial \Lambda_\varepsilon}(\partial \Omega_\varepsilon)$, for $r\ge 0$, the standard Sobolev spaces which are closed subspaces of $H^r(\Omega_\varepsilon)$ and $H^r(\partial \Omega_\varepsilon)$, respectively, and the subscript $\partial \Lambda_\varepsilon$ means that, respectively, traces or functions in $\partial \Omega_\varepsilon$, vanish on this part of the boundary of $\Omega_\varepsilon$, i.e.
$$H^r_{\partial \Lambda_\varepsilon}(\Omega_\varepsilon)=\{v\in H^r(\Omega_\varepsilon):\gamma_0(v)=0  \text{ on } \partial \Lambda_\varepsilon \},$$
and
$$H^r_{\partial \Lambda_\varepsilon}(\partial \Omega_\varepsilon)=\{v\in H^r(\partial \Omega_\varepsilon):v=0  \text{ on } \partial \Lambda_\varepsilon \}.$$
Let us notice that, in fact, we can consider an element of $H^{1/2}(\partial S_\varepsilon)$ as an element of $H^{1/2}_{\partial \Lambda_\varepsilon}(\partial \Omega_\varepsilon)$, and we can consider the given $\psi_\varepsilon^{0}$ as an element of $L^2_{\partial \Lambda_\varepsilon}(\partial \Omega_\varepsilon)$.
%%%%%%%%%%%%%%%%%%

We denote by $H^r_{\partial \Omega}(\widetilde \Omega_\varepsilon)$ and $H^r_{\partial \Omega}(\partial \widetilde\Omega_\varepsilon)$, for $r\ge 0$, the standard Sobolev spaces which are closed subspaces of $H^r(\widetilde \Omega_\varepsilon)$ and $H^r(\partial \widetilde \Omega_\varepsilon)$, respectively, and the subscript $\partial \Omega$ means that, respectively, traces or functions in $\partial \widetilde \Omega_\varepsilon$, vanish on this part of the boundary of $\widetilde \Omega_\varepsilon$, i.e.
$$H^r_{\partial \Omega}(\widetilde\Omega_\varepsilon)=\{\tilde v\in H^r(\widetilde \Omega_\varepsilon):\gamma_0(\tilde v)=0  \text{ on } \partial \Omega\},$$
and
$$H^r_{\partial \Omega}(\partial \widetilde \Omega_\varepsilon)=\{\tilde v\in H^r(\partial \widetilde \Omega_\varepsilon):\tilde v=0  \text{ on } \partial \Omega \}.$$

%%%%%%%%%%%%%%%%%%%%%%%%%%
Let us consider the space
$$H:= L^{2}(
\widetilde \Omega_\varepsilon) \times L_{\partial \Omega}^{2}( \partial\widetilde \Omega_\varepsilon)
\text{,}
$$with the natural inner product $ ((
\tilde v,\tilde \phi), ( \tilde w,\tilde \varphi))_{H}=(\tilde v,\tilde w)_{\widetilde \Omega_\varepsilon}+\varepsilon
(\tilde \phi,\tilde \varphi)_{\partial F_\varepsilon},$ which in particular
induces the norm $|(\cdot,\cdot)|_{H}$ given by
$$|(
\tilde v,\tilde \phi)|^2_{H}=|\tilde v|_{\widetilde \Omega_\varepsilon}^2+\varepsilon|\tilde \phi|^2_{\partial F_\varepsilon},\quad(\tilde v,\tilde \phi)\in
H.$$
Let us also consider the space
$$V_1:=\left\{ \left( \tilde v,\gamma_{0}(\tilde v)\right) :\tilde v\in H_{\partial \Omega}^{1}(
\widetilde \Omega_\varepsilon) \right\}.$$ We note that $V_1$ is a closed
vector subspace of $H_{\partial \Omega}^{1}(\widetilde \Omega_\varepsilon) \times
H_{\partial \Omega}^{1/2}(\partial\widetilde \Omega_\varepsilon),$ and therefore, with
the norm $\|(\cdot,\cdot)\|_{V_1}$ given by
\[
\left\Vert \left( \tilde v,\gamma_{0}(\tilde v)\right) \right\Vert^2
_{V_1}=\left\Vert
\tilde v\right\Vert _{\widetilde \Omega_\varepsilon}^2+\left\Vert \gamma_{0}%
(\tilde v)\right\Vert^2_{\partial F_\varepsilon}, \quad\left(
\tilde v,\gamma_{0}(\tilde v)\right) \in V_1,
\]
$V_1$ is a Hilbert space.

For a vectorial function $v=(v',v_3)$ and a scalar function $w$, we introduce the operators: $\nabla_\varepsilon$ and ${\rm div}_{\varepsilon}$, by
\begin{equation}\label{definition_divergencia_epsilon}
\nabla_{\varepsilon}w=(\nabla_{x^{\prime}}w,\frac{1}{\varepsilon}\partial_{y_3}w)^t, \quad \quad {\rm div}_{\varepsilon}v={\rm div}_{x^{\prime}}v^{\prime}+\frac{1}{\varepsilon}\partial_{y_3}v_3.
\end{equation}
We denote by $|\mathcal{O}|$ the Lebesgue measure of $|\mathcal{O}|$ (3-dimensional if $\mathcal{O}$ is a 3-dimensional open set, 2-dimensional of $\mathcal{O}$ is a curve).

Finally, we denote by $C$ a generic positive constant, independent of $\varepsilon$, which can change from line to line.

\section{Existence and uniqueness of solution of the problem (\ref{PDE})}\label{S3}
We state in this section a result on the existence and uniqueness of solution of problem (\ref{PDE}). 

\begin{definition}\label{definition_weakSolution} A weak solution of (\ref{PDE}) is a pair of functions $(u_\varepsilon,\psi_\varepsilon)$, satisfying
\begin{equation}\label{weak0}
 u_\varepsilon\in
C([0,T];L^2(\Omega_\varepsilon)),\quad \psi_\varepsilon\in
C([0,T];L_{\partial \Lambda_\varepsilon}^2(\partial\Omega_\varepsilon)),\quad\hbox{
for all $T>0$,}
\end{equation}
\begin{equation}\label{weak1}
 u_\varepsilon\in L^2(0,T;H_{\partial \Lambda_\varepsilon}^1(\Omega_\varepsilon)),
 \quad\hbox{
for all $T>0$,}
\end{equation}
\begin{equation}\label{weak2}
\psi_\varepsilon\in L^2(0,T;H_{\partial \Lambda_\varepsilon}^{1/2}(\partial\Omega_\varepsilon)),\quad\hbox{ for all $T>0$,}
\end{equation}
\begin{equation}\label{weak3}
 \gamma_0(u_\varepsilon(t))=\psi_\varepsilon(t),\quad\hbox{ a.e. $t\in (0,T],$}
 \end{equation}
 \begin{equation}\label{weak4}
\left\{
\begin{array}{l}
 \dfrac{d}{dt}(u_\varepsilon(t),v)_{\Omega_\varepsilon}+\varepsilon\,\dfrac{d}{dt}(
\psi_\varepsilon(t),\gamma_{0}(v))_{\partial S_\varepsilon}+(\nabla u_\varepsilon(t),\nabla v)_{
\Omega_\varepsilon}+\kappa\,(u_\varepsilon(t),v)_{\Omega_\varepsilon}\\[2ex]
 =(f(x',t),v)_{\Omega_\varepsilon}
 +\varepsilon\,(g(x',t),\gamma_{0}(v))_{\partial S_\varepsilon}\\[2ex]
 \hbox{in $\mathcal{D}'(0,T)$, for all $v\in H_{\partial \Lambda_\varepsilon}^1(\Omega_\varepsilon),$}
\end{array}
\right.
\end{equation}
\begin{equation}\label{weak5}
 u_\varepsilon(0)=u_\varepsilon^0,\quad and\quad \psi_\varepsilon(0)=\psi_\varepsilon^0.
\end{equation}
\end{definition}
Thanks to \cite[Theorem 3.3]{Anguiano_MJOM}, we have the following result.
\begin{theorem}
\label{Existence_solution_PDE}Under the assumptions (\ref{hyp
0})--(\ref{hyp 0'}), there exists a unique solution
$(u_\varepsilon,\psi_\varepsilon)$
of the problem (\ref{PDE}). Moreover, this solution satisfies the
energy equality
\begin{eqnarray*}
&&\frac{1}{2}\frac{d}{dt}\left(|u_\varepsilon(t)|^2_{\Omega_\varepsilon}+\varepsilon|\psi_\varepsilon(t))|^2_{\partial S_\varepsilon}\right)+|\nabla
u_\varepsilon(t)|^2_{\Omega_\varepsilon} +\kappa\,|u_\varepsilon(t)|^2_{\Omega_\varepsilon} \\\nonumber
&=&(f(x',t),u_\varepsilon(t))_{\Omega_\varepsilon}+\varepsilon\,(g(x',t),\psi_\varepsilon(t))_{\partial S_\varepsilon},\quad\mbox{a.e. $t\in(0,T).$}
\end{eqnarray*}
\end{theorem}
%%%%%%%%%%%%%%%%
\section{Dilatation in the variable $x_3$: definition of $\tilde u_\varepsilon$}\label{S4}
In this section, we use the dilatation in the variable $x_3$, i.e. we use (\ref{dilatacion}) in order to have the functions defined in the open set with fixed height $\widetilde \Omega_\varepsilon$ defined by (\ref{OmegaTilde}). In this sense, we define $\tilde{u}_{\varepsilon}\in L^2(0,T;H_{\partial \Omega}^1(\widetilde \Omega_\varepsilon))$, for all $T>0$, by (\ref{definition_u_tilde}), and using the transformation (\ref{dilatacion}), the system (\ref{PDE}) can be rewritten as
\begin{equation}
\left\{
\begin{array}
[c]{r@{\;}c@{\;}ll}%
\displaystyle\frac{\partial \tilde u_\varepsilon}{\partial t}-{\rm div}_{\varepsilon}(\nabla_\varepsilon \tilde u_\varepsilon)+\kappa\, \tilde u_\varepsilon &
= &
f(x',t)\quad & \text{\ in }\;\widetilde\Omega_\varepsilon\times(0,T) ,\\
\nabla_\varepsilon \tilde u_\varepsilon\cdot \tilde\nu_\varepsilon+\varepsilon\,\displaystyle\frac{\partial \tilde u_\varepsilon}{\partial
t} & = & \varepsilon\,g(x',t) & \text{\ on }%
\;\partial F_\varepsilon\times( 0,T),\\
\tilde u_\varepsilon(x',y_3,0) & = & \tilde u_\varepsilon^{0}(x',y_3), & \text{\ for }\;(x',y_3)\in\widetilde\Omega_\varepsilon,\\
\tilde u_\varepsilon(x',y_3,0) & = & \tilde \psi_\varepsilon^{0}(x',y_3), & \text{\ for
}\;(x',y_3)\in\partial F_\varepsilon,\\
\tilde u_\varepsilon&=& 0, & \text{\ on }%
\;\partial \Omega\times( 0,T),
\end{array}
\right. \label{PDE_dilatado}%
\end{equation}
where ${\rm div}_{\varepsilon}$ and $\nabla_\varepsilon$ are given by (\ref{definition_divergencia_epsilon}), and $\tilde \nu_\varepsilon$ is the outer normal to $\partial F_\varepsilon$.

Taking in (\ref{PDE}) as test function $\tilde v(x',x_3/\varepsilon)$ with $\tilde v\in H_{\partial \Omega}^1(\widetilde \Omega_\varepsilon)$, applying the change of variables (\ref{dilatacion}) and taking into account that $dx=\varepsilon dx'dy_3$ and $d\sigma(x)=\varepsilon d\sigma(x')dy_3$, the variational formulation of problem (\ref{PDE_dilatado}) is then the following one
\begin{equation}\label{weak_dilatado}
\left\{
\begin{array}{l}
 \dfrac{d}{dt}\left(\displaystyle\int_{\widetilde \Omega_\varepsilon}\tilde u_\varepsilon(t) \tilde vdx'dy_3 \right)+\varepsilon\,\dfrac{d}{dt}\left(\displaystyle\int_{\partial F_\varepsilon}\gamma_0(\tilde u_\varepsilon(t))\gamma_{0}(\tilde v)d\sigma(x')dy_3 \right)
 +\displaystyle\int_{\widetilde \Omega_\varepsilon}\nabla_\varepsilon \tilde u_\varepsilon(t)\cdot\nabla_\varepsilon \tilde vdx'dy_3\\[2ex]
 +\kappa\,\displaystyle\int_{\widetilde \Omega_\varepsilon}\tilde u_\varepsilon(t)\tilde v dx'dy_3
 =\displaystyle \int_{\widetilde \Omega_\varepsilon}f(x',t)\tilde vdx'dy_3+\varepsilon\,\displaystyle \int_{\partial F_\varepsilon}g(x',t)\gamma_{0}(\tilde v)d\sigma(x')dy_3\\[2ex]
 \hbox{in $\mathcal{D}'(0,T)$, for all $\tilde v\in H_{\partial \Omega}^1(\widetilde \Omega_\varepsilon).$}
\end{array}
\right.
\end{equation}
Moreover, $\tilde u_\varepsilon$ satisfies the
energy equality
\begin{eqnarray}\label{energyequality_thin}
&&\frac{1}{2}\frac{d}{dt}\left(|(\tilde u_\varepsilon(t),\gamma_0(\tilde u_\varepsilon(t)))|^2_{H}\right)+|\nabla_\varepsilon
\tilde u_\varepsilon(t)|^2_{\widetilde \Omega_\varepsilon} +\kappa\,|\tilde u_\varepsilon(t)|^2_{\widetilde \Omega_\varepsilon} \\
&=&(f(x',t),\tilde u_\varepsilon(t))_{\widetilde \Omega_\varepsilon}+\varepsilon\,(g(x',t),\gamma_0(\tilde u_\varepsilon(t)))_{\partial F_\varepsilon},\quad\mbox{a.e. $t\in(0,T).$}\nonumber
\end{eqnarray}
Considering in (\ref{Initial_condition}) the change of variables given in (\ref{dilatacion}), and taking into account that $dx=\varepsilon dx'dy_3$ and $d\sigma(x)=\varepsilon d\sigma(x')dy_3$, we can deduce that
\begin{equation}\label{Initial_condition_tilde}
|\tilde u_\varepsilon^0|^2_{\widetilde \Omega_\varepsilon}+\varepsilon |\tilde \psi_\varepsilon^0|^2_{\partial F_\varepsilon}\leq C.
\end{equation}

\section{{\it A priori} estimates for $\tilde u_\varepsilon$}\label{S5}
Let us obtain some {\it a priori} estimates for $\tilde u_\varepsilon$.
\begin{lemma}\label{estimates1}
Under the assumptions (\ref{hyp
0})--(\ref{hyp 0'}), there exists a positive constant $C$ independent of $\varepsilon$, such that the solution $\tilde u_\varepsilon$ of the problem (\ref{PDE_dilatado}) satisfies
\begin{eqnarray}\label{acotacion2}
|\nabla_\varepsilon \tilde u_\varepsilon|_{\widetilde \Omega_\varepsilon,T}\leq C, \quad |\tilde u_\varepsilon|_{\widetilde \Omega_\varepsilon,T}\leq C,\quad
\end{eqnarray}
\begin{eqnarray}\label{acotacion_gronwall}
|\tilde u_\varepsilon(t)|_{\widetilde \Omega_\varepsilon}\leq C,\quad \sqrt{\varepsilon}|\gamma_0(\tilde u_\varepsilon(t))|_{\partial F_\varepsilon} \leq C,
\end{eqnarray}
for all $t\in (0,T)$.
\end{lemma}
\begin{proof}
By (\ref{energyequality_thin}), we have
\begin{eqnarray*}
\frac{d}{dt}\left(|(\tilde u_\varepsilon(t),\gamma_0(\tilde u_\varepsilon(t)))|^2_{H}\right)+2|\nabla_\varepsilon
\tilde u_\varepsilon(t)|^2_{\widetilde \Omega_\varepsilon} +2\kappa\,|\tilde u_\varepsilon(t)|^2_{\widetilde \Omega_\varepsilon} 
\leq|(\tilde u_\varepsilon(t),\gamma_0(\tilde u_\varepsilon(t)))|^2_{H}+|f(t)|^2_{\omega_\varepsilon}+\varepsilon |g(t)|^2_{\partial F'_\varepsilon}.
\end{eqnarray*}
Integrating between $0$ and $t$, we obtain
\begin{eqnarray}\label{equality_Gronwall}
&&|(\tilde u_\varepsilon(t),\gamma_0(\tilde u_\varepsilon(t)))|^2_{H}+2{\rm min}\{1,\kappa\}\int_0^t\left(|\nabla_\varepsilon
\tilde u_\varepsilon(s)|^2_{\widetilde \Omega_\varepsilon} +|\tilde u_\varepsilon(s)|^2_{\widetilde \Omega_\varepsilon} \right)ds\\
&\leq&|(\tilde u_\varepsilon^0,\gamma_0(\tilde u_\varepsilon^0))|^2_H+\int_0^t|(\tilde u_\varepsilon(s),\gamma_0(\tilde u_\varepsilon(s)))|^2_{H}ds+\int_0^T\left(|f(s)|^2_{ \omega_\varepsilon}+\varepsilon |g(s)|^2_{\partial F'_\varepsilon}\right)ds.\nonumber
\end{eqnarray}
By \cite[Lemma 4.1]{Anguiano_MJOM} with $p=2$, we can deduce
\begin{equation*}
\varepsilon\,|g(s)|^2_{\partial F'_\varepsilon}\leq C\left(|g(s)|^2_{\omega_\varepsilon}+\varepsilon^2 |\nabla g(s)|^2_{\omega_\varepsilon}\right)\leq C ||g(s)||^2_{\omega_\varepsilon},
\end{equation*}
which together with (\ref{hyp 0'}) gives
\begin{equation}\label{trace1}
\int_0^T\left(|f(s)|^2_{\omega_\varepsilon}+\varepsilon\,|g(s)|^2_{\partial F'_\varepsilon}
\right)ds\leq C.
\end{equation}
Taking into account (\ref{Initial_condition_tilde}) and (\ref{trace1}) in (\ref{equality_Gronwall}) and applying Gronwall Lemma, in particular we obtain (\ref{acotacion_gronwall}). Finally, taking into account (\ref{Initial_condition_tilde}), (\ref{acotacion_gronwall}) and (\ref{trace1}) in (\ref{equality_Gronwall}), we get (\ref{acotacion2}).
\end{proof}

Now, if we want to take the inner product in (\ref{PDE_dilatado}) with $\tilde u'_\varepsilon$, we need that $u'_\varepsilon \in L^2(0,T;H_{\partial \Omega}^1(\widetilde \Omega_\varepsilon))$. However, we do not have it for our weak solution. Therefore, we use the Galerkin method in order to prove, rigorously, new {\it a priori} estimates for $\tilde u_\varepsilon$.

On the space $V_1$ we define a continuous symmetric linear
operator $A_1:V_1\rightarrow V_1^{\prime}$, given by
\begin{equation}\label{def_A1}
\langle A_1(( \tilde v,\gamma_{0}(\tilde v))) ,( \tilde w,\gamma _{0}(\tilde w))\rangle
=(\nabla_\varepsilon \tilde v,\nabla_\varepsilon \tilde w)_{\widetilde \Omega_\varepsilon}+\kappa\,(\tilde v,\tilde w)_{\widetilde \Omega_\varepsilon}\text{, \
}
\end{equation}
for all $\tilde v,\tilde w\in H_{\partial \Omega}^{1}( \widetilde \Omega_\varepsilon)$.

We observe that $A_1$ is coercive. In fact, taking into account that $\varepsilon<1$, we have
\begin{eqnarray}\label{Coercitivity}
\left\langle A_1\left( \left( \tilde v,\gamma_{0}(\tilde v)\right)
,\left( \tilde v,\gamma _{0}(\tilde v)\right) \right) \right\rangle &=&|\nabla_{x'}\tilde v|^2_{\widetilde \Omega_\varepsilon}+{1\over \varepsilon^2}|\partial_{y_3}\tilde v|^2_{\widetilde \Omega_\varepsilon}+\kappa\,|\tilde v|^2_{\widetilde \Omega_\varepsilon}\\
&>& |\nabla_{x'}\tilde v|^2_{\widetilde \Omega_\varepsilon}+|\partial_{y_3}\tilde v|^2_{\widetilde \Omega_\varepsilon}+\kappa\,|\tilde v|^2_{\widetilde \Omega_\varepsilon}\nonumber
\\
&
\ge&\min\left\{ 1,\kappa\right\}
\left\Vert \tilde v\right\Vert _{\widetilde \Omega_\varepsilon}^{2}\nonumber
\\
& = &\frac{1}{1+\|\gamma_0\|^2}\min\left\{ 1,\kappa\right\}
\left\Vert \tilde v\right\Vert
_{ \widetilde \Omega_\varepsilon }^{2}\nonumber\\
&& +\frac{\|\gamma_0\|^2}{1+\left\Vert
\gamma_{0}\right\Vert ^{2}}\min\left\{
1,\kappa\right\} \left\Vert \tilde v\right\Vert _{\widetilde \Omega_\varepsilon}^{2}\nonumber\\
& \geq&\frac{1}{1+\|\gamma_0\|^2}\min\left\{
1,\kappa\right\}
 \left\Vert \left( \tilde v,\gamma_{0}(\tilde v)\right) \right\Vert _{V_1}%
^{2}\text{,}\nonumber
\end{eqnarray}for all $\tilde v\in H_{\partial \Omega}^1(\widetilde\Omega_\varepsilon)$.

Let us observe that the space $H^{1}_{\partial \Omega}(\widetilde \Omega_\varepsilon)\times H_{\partial \Omega}^{1/2}(\partial\widetilde \Omega_\varepsilon)$ is compactly imbedded in $H$, and
therefore, for the symmetric and coercive linear continuous operator $A_{1}:V_{1}\rightarrow V_{1}^{\prime}$, where $A_1$ is given by (\ref{def_A1}), there exists a non-decreasing sequence $0<\lambda_{1}\leq\lambda_{2}\leq\ldots$ of eigenvalues associated to the operator
$A_{1}$ with $\lim_{j\rightarrow\infty }\lambda_{j}=\infty,$ and there exists a Hilbert basis of $H$, $\{(w_{j},\gamma_{0}(w_{j})) :j\geq1\}$$\subset D(A_1)$, with $span\{(w_{j},\gamma_{0}(w_{j})):j\geq1\} $ densely embedded in $V_{1}$, such that
\[
A_{1}((w_{j},\gamma_{0}(w_{j})))=\lambda_{j}(w_{j},\gamma_{0}(w_{j}))\quad\forall j\geq1.
\]
Taking into account the above facts, we denote by $$(\tilde u_{\varepsilon,m}(t),\gamma_{0}(\tilde u_{\varepsilon,m}(t)))=(\tilde u_{\varepsilon,m}(t;0,\tilde u_\varepsilon^0,\tilde \psi_\varepsilon^0),\gamma_{0}(\tilde u_{\varepsilon,m}(t;0,\tilde u_\varepsilon^0,\tilde \psi_\varepsilon^0)))$$ the Galerkin approximation of the solution
$(\tilde u_\varepsilon(t;0,\tilde u_\varepsilon^0,\tilde \psi_\varepsilon^0),\gamma_{0}(\tilde u_\varepsilon(t;0,\tilde u_\varepsilon^0,\tilde \psi_\varepsilon^0)))$ to (\ref{PDE_dilatado}) for each integer $m\geq1$, which is given by
\begin{equation}
(\tilde u_{\varepsilon,m}(t),\gamma_{0}(\tilde u_{\varepsilon,m}(t)))=\sum_{j=1}^{m}\delta_{\varepsilon mj}(t)(w_{j},\gamma_{0}(w_{j})),\label{Galerkin1}
\end{equation}
and is the solution of
\begin{eqnarray}
\nonumber &&\dfrac{d}{dt}((\tilde u_{\varepsilon,m}(t),\gamma_{0}(\tilde u_{\varepsilon,m}(t))),(w_{j},\gamma_{0}(w_{j})))_{H}\\ \nonumber
&&+\left\langle A_{1}((\tilde u_{\varepsilon,m}(t),\gamma_{0}(\tilde u_{\varepsilon,m}(t)))),(w_{j},\gamma_{0}(w_{j}))\right\rangle\\ 
&=&( f(t),w_{j})_{\widetilde \Omega_\varepsilon}+\varepsilon( g(t),\gamma_{0}(w_{j}))_{\partial F_\varepsilon},\quad j=1,\ldots ,m,\label{7}
\end{eqnarray}
with initial data
\begin{equation}
\label{7'}
(\tilde u_{\varepsilon,m}(0),\gamma_{0}(\tilde u_{\varepsilon,m}(0)))=(\tilde u_{\varepsilon,m}^{0},\gamma_{0}(\tilde u_{\varepsilon,m}^{0})),
\end{equation}
where
\[
\delta_{\varepsilon mj}(t)=(\tilde u_{\varepsilon,m}(t),w_{j})_{\widetilde\Omega_\varepsilon}+( \gamma_{0}(\tilde u_{\varepsilon,m}(t)),\gamma_{0}(w_{j}))_{\partial F_\varepsilon},
\]
and $(\tilde u_{\varepsilon,m}^{0},\gamma_0(\tilde u_{\varepsilon,m}^{0}))\in span\{(w_j,\gamma_0(w_j)): j=1,\ldots ,m\}$ converge (when
$m\to\infty$) to $(\tilde u_\varepsilon^0,\tilde \psi_\varepsilon^0)$ in a suitable sense which will be specified below.
\begin{lemma}\label{Lemmaestimates}
Suppose the assumptions (\ref{hyp
0})--(\ref{hyp 0'}). Then, for any initial condition $(u_\varepsilon^0,\psi_\varepsilon^0)\in H^1(\Omega_\varepsilon)\times H^{1/2}(\partial S_\varepsilon)$ of the problem (\ref{PDE}) such that satisfies (\ref{Initial_condition_2}), there exists a positive constant $C$ independent of $\varepsilon$, such that the solution $\tilde u_\varepsilon$ of the problem (\ref{PDE_dilatado}) satisfies
\begin{eqnarray}\label{acotacion5}
 \sup_{t\in [0,T]}\left\Vert \tilde u_\varepsilon(t)\right\Vert _{\widetilde\Omega_\varepsilon}\leq C,
\end{eqnarray}
\begin{eqnarray}\label{acotacion3}
|\tilde u'_\varepsilon|_{\widetilde \Omega_\varepsilon,T} \leq C, \quad \sqrt{\varepsilon}|\gamma_0(\tilde u'_\varepsilon)|_{\partial F_\varepsilon,T}\leq C.
\end{eqnarray}
\end{lemma}
\begin{proof}
First, considering in (\ref{Initial_condition_2}) the change of variables given in (\ref{dilatacion}), and taking into account that $dx=\varepsilon dx'dy_3$ and $\partial_{y_3}=\varepsilon \partial_{x_3}$, we can deduce that $(\tilde u_\varepsilon^0,\tilde \psi_\varepsilon^0)\in V_1$ such that
\begin{equation}\label{Initial_condition_2_tilde}
|\nabla_\varepsilon \tilde u_\varepsilon^0|_{\widetilde \Omega_\varepsilon}\leq C.
\end{equation}
For all $m\geq1$, there exists $(\tilde u_{\varepsilon,m}^{0},\gamma_{0}(\tilde u_{\varepsilon,m}^{0}))\in span\{(w_{j},\gamma_{0}(w_{j})):1\leq
j\leq m\} $, such that the sequence $\{(\tilde u_{\varepsilon,m}^{0},\gamma_{0}(\tilde u_{\varepsilon,m}^{0}))\}$ converges to $(\tilde u_\varepsilon^0,\tilde \psi_\varepsilon^0)$ in $V_1$. Then, taking into account (\ref{Initial_condition_tilde}) and (\ref{Initial_condition_2_tilde}), we know that there exists a positive constant $C$ such that
\begin{equation}\label{Initial_condition_2_Galerkin}
|\nabla_\varepsilon \tilde u_{\varepsilon,m}^0|^2_{\widetilde \Omega_\varepsilon}+|\tilde u_{\varepsilon,m}^0|^2_{\widetilde \Omega_\varepsilon}\leq C.
\end{equation}
For each integer $m\geq1$, we consider the sequence $\{(u_{\varepsilon,m}(t),\gamma_{0}(u_{\varepsilon,m}(t)))\}$ defined by
(\ref{Galerkin1})-(\ref{7'}) with these initial data. Multiplying by the derivative $\delta'_{\varepsilon mj}$ in (\ref{7}), and summing from $j=1$ to $m$, we obtain
\begin{eqnarray*}
&& |(\tilde u_{\varepsilon,m}^{\prime}(t),\gamma_{0}(\tilde u_{\varepsilon,m}^{\prime}(t)))|^2_{H}\\
&&+\frac{1}{2}\frac{d}{dt}(\left\langle A_{1}((\tilde u_{\varepsilon,m}(t),\gamma_{0}(\tilde u_{\varepsilon,m}(t)))),(\tilde u_{\varepsilon,m}(t),\gamma_{0}(\tilde u_{\varepsilon,m}(t)))\right\rangle)\\
&  =&(f(t),\tilde u_{\varepsilon,m}^{\prime}(t))_{\widetilde \Omega_\varepsilon}+\varepsilon(g(t),\gamma_{0}(\tilde u_{\varepsilon,m}^{\prime}(t)))_{\partial F_\varepsilon}.
\end{eqnarray*}
Then, we deduce
\begin{eqnarray*}\nonumber
&&|(\tilde u_{\varepsilon,m}^{\prime}(t),\gamma_{0}(\tilde u_{\varepsilon,m}^{\prime}(t)))|^2_{H}\\
&&+\frac{1}{2}\frac{d}{dt}(\left\langle A_{1}((\tilde u_{\varepsilon,m}(t),\gamma_{0}(\tilde u_{\varepsilon,m}(t)))),(\tilde u_{\varepsilon,m}(t),\gamma_{0}(\tilde u_{\varepsilon,m}(t)))\right\rangle)\\
&  \leq& \frac{1}{2}|f(t)|_{\omega_\varepsilon}^{2}+\frac{1}{2}|\tilde u_{\varepsilon,m}^{\prime}(t)|_{\widetilde \Omega_\varepsilon}^{2}+\varepsilon\frac{1}{2}|g(t)|_{\partial F'_\varepsilon}^{2} +\varepsilon\frac{1}{2}|\gamma_{0}(\tilde u_{\varepsilon,m}^{\prime}(t))|_{\partial F_\varepsilon}^{2}.\label{acotacion_necesaria}
\end{eqnarray*}
Integrating now between $0$ and $t$, taking into account the definition of $A_1$ and (\ref{Coercitivity}), we obtain that
\begin{eqnarray}\label{last_estimate}
&& \int_{0}^{t}|(\tilde u_{\varepsilon,m}^{\prime}(s),\gamma_{0}(\tilde u_{\varepsilon,m}^{\prime}(s)))|^2_{H}ds+ \frac {\min\{1,\kappa\}}{1+\|\gamma_{0}\|^{2}}\|(\tilde u_{\varepsilon,m}(t),\gamma_{0}(\tilde u_{\varepsilon,m}(t)))\|_{V_{1}}^{2}\nonumber\\
&\leq&\max\{1,\kappa\}\left(|\nabla_\varepsilon \tilde u_{\varepsilon,m}^{0}|^2_{\widetilde \Omega_\varepsilon}+|\tilde u_{\varepsilon,m}^{0}|^2_{\widetilde \Omega_\varepsilon}\right)
+\int_{0}^{T}(|f(s)|_{\omega_\varepsilon}^{2}+\varepsilon |g(s)|_{\partial F'_\varepsilon}^{2})ds,
\end{eqnarray}
for all $t\in (0,T)$.

Taking into account (\ref{trace1}) and (\ref{Initial_condition_2_Galerkin}) in (\ref{last_estimate}), we have proved that the sequence $\{(\tilde u_{\varepsilon,m},\gamma_0(\tilde u_{\varepsilon,m}))\}$ is bounded in $C([0,T];V_1),$  and $\{(\tilde u_{\varepsilon,m}',\gamma_0(\tilde u_{\varepsilon,m}'))\}$ is bounded in $L^2(0, T;H),$
for all $T>0$.

If we work with the truncated Galerkin equations (\ref{Galerkin1})-(\ref{7'}) instead of the full PDE, using (\ref{Coercitivity}), we note that the calculations of the proof of Lemma \ref{estimates1} can be followed identically to show that $\{(\tilde u_{\varepsilon,m},\gamma_0(\tilde u_{\varepsilon,m}))\}$ is bounded in $L^2(0,T;V_1),$
for all $T>0$.

Moreover, taking into account the uniqueness of solution to (\ref{PDE_dilatado}) and using Aubin-Lions compactness lemma (e.g., cf. Lions \cite{Lions}), it is not difficult to conclude that 
the sequence $\{( \tilde u_{\varepsilon,m},\gamma_{0} (\tilde u_{\varepsilon, m}))\}$ converges weakly in $
L^{2}(0,T;V_1)$ to the solution
$(\tilde u_\varepsilon,\gamma_{0}(\tilde u_\varepsilon))$ to (\ref{PDE_dilatado}). Since the inclusion $H^1(\widetilde \Omega_\varepsilon)\subset L^2(\widetilde \Omega_\varepsilon)$ is compact and $\tilde u_\varepsilon\in C([0,T];L^2(\widetilde \Omega_\varepsilon))$, it follows using \cite[Lemma 11.2]{Robinson} that the estimate (\ref{acotacion5}) is proved.

On the other hand, the sequence $\{( \tilde u'_{\varepsilon,m},\gamma_{0} (\tilde u'_{\varepsilon, m}))\}$ converges weakly in $
L^{2}(0,T;H)$ to $(\tilde u'_\varepsilon,\gamma_{0}(\tilde u'_\varepsilon))$, for all $T>0$, and using the lower-semicontinuity of the norm and (\ref{last_estimate}), we get
\begin{eqnarray*}
|\tilde u'_\varepsilon|^2_{\widetilde \Omega_\varepsilon,T}+\varepsilon |\gamma_0(\tilde u'_{\varepsilon})|^2_{\partial F_\varepsilon,T}&\leq& \liminf_{m\to \infty}\left(|\tilde u'_{\varepsilon,m}|^2_{\widetilde \Omega_\varepsilon,T}+\varepsilon |\gamma_0(\tilde u'_{\varepsilon,m})|^2_{\partial F_\varepsilon,T}\right)\\
&\leq& C\liminf_{m\to \infty}\left(|\nabla_\varepsilon \tilde u_{\varepsilon,m}^{0}|^2_{\widetilde \Omega_\varepsilon}+|\tilde u_{\varepsilon,m}^{0}|^2_{\widetilde \Omega_\varepsilon}+1\right) \\
&=&C\left(|\nabla_\varepsilon \tilde u_{\varepsilon}^{0}|^2_{\widetilde \Omega_\varepsilon}+|\tilde u_{\varepsilon}^{0}|^2_{\widetilde \Omega_\varepsilon}+1\right),
\end{eqnarray*}
which, jointly with (\ref{Initial_condition_tilde}) and (\ref{Initial_condition_2_tilde}), implies (\ref{acotacion3}).
\end{proof}

\subsection{The extension of $\tilde u_\varepsilon$ to the whole $\Omega \times (0,T)$: definition of $\widetilde \varPi_\varepsilon \tilde u_\varepsilon$}
Since the solution $\tilde u_\varepsilon$ of the problem (\ref{PDE_dilatado}) is defined only in $\widetilde \Omega_\varepsilon\times (0,T)$, we need to extend it to the whole $\Omega\times (0,T)$. For finding a suitable extension into all $\Omega\times (0,T)$, we shall use the following result.
\begin{corollary}[Corollary 3 in \cite{Anguiano_SuarezGrau}]\label{corollary_para_Poincare}
There exists an extension operator $\widetilde \varPi_\varepsilon\in \mathcal{L}(H^1_{\partial \Omega}(\widetilde \Omega_\varepsilon);H_0^1(\Omega))$ and a positive constant $C$, independent of $\varepsilon$, such that 
\begin{equation*}
\widetilde \varPi_\varepsilon \tilde \varphi(x',y_3)=\tilde \varphi (x',y_3), \quad \text{ if \ } (x',y_3)\in \widetilde \Omega_\varepsilon,
\end{equation*}
\begin{equation*}
|\nabla_\varepsilon \widetilde \varPi_\varepsilon \tilde \varphi |_{\Omega}\leq C|\nabla_\varepsilon  \tilde \varphi |_{\widetilde \Omega_\varepsilon}, \quad \forall\, \tilde \varphi \in H^1_{\partial \Omega}(\widetilde \Omega_\varepsilon).
\end{equation*}
\end{corollary}
Let us obtain some a priori estimates for the extension of $\tilde u_\varepsilon$ to the whole $\Omega\times (0,T)$. Using Corollary \ref{corollary_para_Poincare} together with Lemmas \ref{estimates1}-\ref{Lemmaestimates}, we obtain the following result.
\begin{corollary}\label{estimates_extension}
Assume the assumptions in Lemma \ref{Lemmaestimates}. Then, there exists an extension $\widetilde \varPi_\varepsilon \tilde u_\varepsilon$ of the solution $\tilde u_\varepsilon$ of the problem (\ref{PDE_dilatado}) into $\Omega\times (0,T)$, such that
\begin{eqnarray}\label{acotacion1_extension}
|\nabla_\varepsilon \widetilde \varPi_\varepsilon\tilde u_\varepsilon|_{ \Omega,T}\leq C, \quad |\widetilde \varPi_\varepsilon\tilde u_\varepsilon|_{ \Omega,T}\leq C,\quad
\end{eqnarray}
\begin{equation}\label{acotacion2_extension}
\sup_{t\in [0,T]}\left\Vert  \widetilde \varPi_\varepsilon\tilde u_\varepsilon(t)\right\Vert _{\Omega}\leq C, 
\end{equation}
\begin{eqnarray}\label{acotacion3_extension}
|\widetilde \varPi_\varepsilon\tilde u'_\varepsilon|_{ \Omega,T} \leq C,
\end{eqnarray}
where the constant $C$ does not depend on $\varepsilon$.
\end{corollary}

\section{A compactness result}\label{S6}
In this section, we obtain some compactness results about the behavior of the sequence $\widetilde \varPi_\varepsilon \tilde u_\varepsilon$ satisfying the {\it a priori} estimates given in Corollary \ref{estimates_extension}.

Due to the periodicity of the domain $\widetilde \Omega_\varepsilon$, one has, for $\varepsilon \to 0$, that 
\begin{equation}\label{convergence_chi}
\chi_{\widetilde\Omega_\varepsilon}\stackrel{\tt
*}\rightharpoonup {|Y'_f|\over |Y'|} \quad \textrm{weakly-star in}\
L^\infty(\Omega),
\end{equation}
where the limit is the proportion of the material in the cell $Y'$.

Set
$$
\xi_\varepsilon=\nabla_{x'} \tilde u_\varepsilon\quad \text{ in }\widetilde\Omega_\varepsilon\times (0,T),
$$
which, by the first estimate in (\ref{acotacion2}), satisfies
\begin{equation}\label{estimate_xi_epsilon}
|\xi_\varepsilon|_{\widetilde \Omega_\varepsilon,T}\leq C.
\end{equation}
Let us denote by $\tilde \xi_\varepsilon$ its extension with zero to the whole of $\Omega \times (0,T)$, i.e.
\begin{equation}\label{definition_tildexi}
\tilde \xi_\varepsilon=\left\{
\begin{array}{l}
\xi_\varepsilon \quad \text{in }\widetilde \Omega_\varepsilon \times (0,T),\\
 0 \quad \text{in }(\Omega\setminus \overline{\widetilde\Omega_\varepsilon})\times (0,T).
 \end{array}\right.
\end{equation}

\begin{proposition}\label{Propo_convergence}
Under the assumptions in Lemma \ref{Lemmaestimates}, there exists a function $u\in L^2(0,T;H_0^1(\omega))$, where $u$ is independent of $y_3$, ($u$ will be the unique solution of the limit system (\ref{limit_problem})) and a function $\xi\in L^2(0,T;(L^2(\Omega))^3)$ such that for all $T>0,$
\begin{eqnarray}
\label{continuity1} 
&\widetilde \varPi_\varepsilon \tilde{u}_\varepsilon(t)\rightharpoonup
 u(t) &\textrm{weakly in}\ H_0^1(\Omega),\quad \forall t\in[0,T],
\\
\label{converge_initial_data}
&\widetilde \varPi_\varepsilon \tilde{u}_{\varepsilon}(t)\rightarrow
u(t)&\quad \text{strongly in }L^2(\Omega),\quad \forall t\in[0,T],
\\
\label{converge_gradiente}
&\tilde \xi_\varepsilon\rightharpoonup \xi& \quad \text{weakly in} \quad L^2(0,T;(L^2(\Omega))^3),
\end{eqnarray}
where $\tilde \xi_\varepsilon$ is given by (\ref{definition_tildexi}).
\end{proposition}

\begin{proof}
The estimates (\ref{acotacion1_extension}) read
\begin{eqnarray}\label{acotacion1_extension_separadas}
|\nabla_{x'} \widetilde \varPi_\varepsilon\tilde u_\varepsilon|_{ \Omega,T}\leq C, \quad |\partial_{y_3} \widetilde \varPi_\varepsilon\tilde u_\varepsilon|_{ \Omega,T}\leq C\,\varepsilon\quad |\widetilde \varPi_\varepsilon\tilde u_\varepsilon|_{ \Omega,T}\leq C.
\end{eqnarray}
From the first and the last estimates in (\ref{acotacion1_extension_separadas}), we can deduce that
$\{\widetilde \varPi_\varepsilon\tilde u_\varepsilon\}$ is bounded in $L^2(0,T;L^2(0,1;H^1(\omega)))$, for all $T>0.$ Let us fix $T>0$. Then, there exist a subsequence $\{\widetilde \varPi_\varepsilon\tilde u_{\varepsilon'}\}\subset \{\widetilde \varPi_\varepsilon\tilde u_\varepsilon\}$ and a function $u\in L^2(0,T;L^2(0,1;H^1(\omega)))$ such that
\begin{eqnarray}
\label{compacidad1_prueba}
\widetilde \varPi_\varepsilon\tilde{u}_{\varepsilon'}\rightharpoonup
u &\textrm{weakly in}\ L^2(0,T;L^2(0,1;H^1(\omega))),
\end{eqnarray}
which implies
\begin{eqnarray}
\label{compacidad2_prueba}
\partial_{y_3}\widetilde \varPi_\varepsilon\tilde{u}_{\varepsilon'}\rightharpoonup
\partial_{y_3}u &\textrm{weakly in}\ L^2(0,T;H^{-1}(0,1;H^1(\omega))).
\end{eqnarray}
By the second estimate in (\ref{acotacion1_extension_separadas}), we can deduce that $\partial_{y_3} \widetilde \varPi_\varepsilon\tilde u_{\varepsilon'}$ tends to zero in $L^2(0,T;L^2(\Omega))$. From (\ref{compacidad2_prueba}) and the uniqueness of the limit, we have that $\partial_{y_3}u=0$, which implies that $u$ does not depend on $y_3$.

On the other hand, by the estimate (\ref{acotacion3_extension}), we see that the sequence
$\{\widetilde \varPi_\varepsilon\tilde u'_\varepsilon\}$ is bounded in $L^2(0,T;L^2(\Omega))$, for all $T>0.$
Then, we have that $\widetilde \varPi_\varepsilon\tilde u_\varepsilon(t):[0,T]\longrightarrow L^2(\Omega)$ is an equicontinuous family of functions.

By the estimate (\ref{acotacion2_extension}), for each $t\in[0,T]$, we have that $\{\widetilde \varPi_\varepsilon\tilde u_\varepsilon(t)\}$ is bounded in $H_0^1(\Omega)$, so that the compact embedding $H_0^1(\Omega)\subset L^2(\Omega)$, implies that it is precompact in $L^2(\Omega)$. Then, applying the Ascoli-Arzelà Theorem, we deduce that $\{\widetilde \varPi_\varepsilon\tilde u_\varepsilon(t)\}$ is a precompact sequence in $C([0,T];L^2(\Omega))$. Hence, since by (\ref{compacidad1_prueba}) we obtain
\begin{eqnarray*}
\widetilde \varPi_\varepsilon\tilde{u}_{\varepsilon'}\rightharpoonup
u &\textrm{weakly in}\ L^2(0,T;L^2(\Omega)),
\end{eqnarray*}
we have
\begin{eqnarray}\label{converge_strongly_p}
 &\widetilde \varPi_\varepsilon\tilde{u}_{\varepsilon'}\rightarrow
u &\textrm{strongly in}\ C([0,T];L^2(\Omega)).
\end{eqnarray}
The boundedness of $\{\widetilde \varPi_\varepsilon\tilde u_\varepsilon(t)\}$ in $H_0^1(\Omega)$ implies then by a standard argument that
\begin{eqnarray*}
 &\widetilde \varPi_\varepsilon\tilde{u}_{\varepsilon'}(t)\rightharpoonup
u (t)&\textrm{weakly in}\ H_0^1(\Omega),\quad \forall t\in[0,T].
\end{eqnarray*}
From (\ref{converge_strongly_p}), in particular, we have
$$\widetilde \varPi_\varepsilon\tilde{u}_{\varepsilon'}(t)\rightarrow
u(t)\quad \text{strongly in }L^2(\Omega), \quad \forall t\in[0,T].$$

Finally, from the estimate (\ref{estimate_xi_epsilon}) and (\ref{definition_tildexi}), we have $|\tilde \xi_\varepsilon|_{\Omega,T}\leq C$, and hence, up a sequence, there exists $\xi\in L^2(0,T,(L^2(\Omega))^3)$ such that $\tilde \xi_{\varepsilon''} \rightharpoonup \xi$ weakly in $L^2(0,T;(L^2(\Omega))^3)$.
  
By the uniqueness of solution of the limit problem (\ref{limit_problem}), we deduce that the above convergences hold for the whole sequence and therefore, by the arbitrariness of $T>0$, all the convergences are satisfied, as we wanted to prove.
\end{proof}

\section{Homogenized model: proof of Theorem \ref{Main}}\label{S7}
In this section, we identify the homogenized model. 

Let $v\in \mathcal{D}(\Omega)$ be a test function in (\ref{weak_dilatado}), with $v$ independent of $y_3$. Then one has
\begin{eqnarray*}
\dfrac{d}{dt}\left(\int_{\Omega}\chi_{\widetilde\Omega_\varepsilon}\widetilde \varPi_\varepsilon\tilde u_\varepsilon(t)vdx'dy_3\right)+\varepsilon\,\dfrac{d}{dt}\left(\int_{\partial F_\varepsilon}
\gamma_{0}(\tilde u_\varepsilon(t))vd\sigma(x')dy_3\right)+\int_{\Omega}\tilde \xi_\varepsilon \cdot\nabla_{x'} vdx'dy_3\\[2ex]
+\kappa \int_{\Omega}\chi_{\widetilde \Omega_\varepsilon}\widetilde \varPi_\varepsilon\tilde u_\varepsilon(t)v dx'dy_3
=\int_{\Omega} \chi_{\widetilde\Omega_\varepsilon} f(t)vdx'dy_3
 +\varepsilon\int_{\partial F_\varepsilon}g(t)vd\sigma(x')dy_3, 
 \end{eqnarray*}
 in $\mathcal{D}'(0,T)$.
 
We consider $\varphi\in C_c^1([0,T])$ such that $\varphi(T)=0$ and $\varphi(0)\ne 0$. Multiplying by $\varphi$ and integrating between $0$ and $T$, we have
 \begin{eqnarray}\label{system1}\nonumber
-\varphi(0)\left(\int_{\Omega}\chi_{\widetilde \Omega_\varepsilon}\widetilde \varPi_\varepsilon\tilde u_\varepsilon(0)vdx'dy_3\right)-\int_0^T\dfrac{d}{dt}\varphi(t)\left(\int_{\Omega}\chi_{\widetilde\Omega_\varepsilon}\widetilde \varPi_\varepsilon\tilde u_\varepsilon(t)vdx'dy_3\right)dt\\[2ex]\nonumber
-\varepsilon\varphi(0)\!\left(\int_{\partial F_\varepsilon}\!\!\!
\gamma_{0}(\tilde u_\varepsilon(0))vd\sigma(x')dy_3\!\right)
\!-\!\varepsilon\!\!\int_0^T \!\!\dfrac{d}{dt}\varphi(t)\left(\int_{\partial F_\varepsilon}\!\!\!
\gamma_{0}(\tilde u_\varepsilon(t))vd\sigma(x')dy_3\!\right)dt\\[2ex]
+\int_0^T\varphi(t)\int_{\Omega}\tilde \xi_\varepsilon\cdot\nabla_{x'} vdx'dy_3dt
+\kappa \int_0^T\varphi(t)\int_{\Omega}\chi_{\widetilde \Omega_\varepsilon}\widetilde \varPi_\varepsilon\tilde u_\varepsilon(t)v dx'dy_3dt \\[2ex]
=\int_0^T \varphi(t)\int_{\Omega}\chi_{\widetilde\Omega_\varepsilon} f(t)vdx'dy_3dt
 +\varepsilon \int_0^T \varphi(t)\int_{\partial F_\varepsilon}g(t)vd\sigma(x')dy_3dt.\nonumber
 \end{eqnarray}

For the sake of clarity, we split the proof in five parts. Firstly, we analyze the integrals on $\Omega$, where we only require to use Proposition \ref{Propo_convergence} and the convergence (\ref{convergence_chi}). Secondly, for the integrals on the boundary of the cylinders we make use of a convergence result based on a technique introduced by Vanninathan \cite{Vanni}. In the third step, we pass to the limit, as $\varepsilon \to 0$ in order to get the limit equation satisfied by $u$. In the fourth step, we identify $\int_0^1\xi dy_3$ making use of the solutions of the cell-problems (\ref{system_eta}) and, finally we prove that $u$ is uniquely determined.

{\bf Step 1}. In this step, we analyze all the integrals on $\Omega$. Passing to the limit, as $\varepsilon\to 0$, in the integrals on $\Omega$ and taking into account that $u$ is independent of $y_3$:

From (\ref{convergence_chi}) and (\ref{converge_initial_data}), we have, for $\varepsilon \to 0$, 
\begin{equation*}
\int_{\Omega}\chi_{\widetilde \Omega_\varepsilon}\widetilde \varPi_\varepsilon\tilde u_\varepsilon(t)v dx'dy_3\to {|Y'_f|\over |Y'|}\int_{\omega}u(t)v dx',
\end{equation*}
which integrating in time and using Lebesgue's Dominated Convergence Theorem, gives
\begin{equation*}
\int_0^T{d \over dt}\varphi(t)\left(\int_{\Omega}\chi_{\widetilde \Omega_\varepsilon}\widetilde \varPi_\varepsilon\tilde u_\varepsilon(t) vdx'dy_3\right)dt\to {|Y'_f|\over |Y'|}\int_0^T{d \over dt}\varphi(t)\left(\int_{\omega}u(t)vdx'\right)dt,
\end{equation*}
and
\begin{equation*}
\kappa\int_0^T\varphi(t)\int_{\Omega}\chi_{\widetilde\Omega_\varepsilon}\widetilde \varPi_\varepsilon\tilde u_\varepsilon(t) vdx'dy_3dt\to \kappa{|Y'_f|\over |Y'|}\int_0^T\varphi(t)\int_{\omega}u(t)vdx'dt.
\end{equation*}
By (\ref{convergence_chi}) and (\ref{converge_initial_data}), we have
\begin{equation*}
\varphi(0)\left(\int_{\Omega}\chi_{\widetilde\Omega_\varepsilon}\widetilde \varPi_\varepsilon\tilde u_\varepsilon(0)vdx'dy_3\right)\to \varphi(0){|Y'_f|\over |Y'|}\int_{\omega} u(0)vdx'.
\end{equation*}
By the assumption (\ref{hyp 0'}), (\ref{convergence_chi}) and using Lebesgue's Dominated Convergence Theorem, we get
$$\int_0^T \varphi(t)\int_{\Omega}\chi_{\widetilde\Omega_\varepsilon} f(x',t)vdx'dy_3dt \to {|Y'_f|\over |Y'|}\int_0^T\varphi(t)\int_{\omega}f(x',t)vdx'dt.$$
On the other hand, using (\ref{converge_gradiente}), we obtain, for $\varepsilon\to 0$
$$\int_0^T\varphi(t)\int_{\Omega}\tilde \xi_\varepsilon\cdot\nabla_{x'} vdx'dy_3dt \to \int_0^T\varphi(t)\int_{\Omega}\xi\cdot\nabla_{x'} vdx'dy_3dt.$$

{\bf Step 2}. Passing to the limit, as $\varepsilon\to 0$, in the surface integrals on the boundary of the cylinders and taking into account that $u$ is independent of $y_3$:

We make use of the technique introduced by Vanninathan \cite{Vanni} for the Steklov problem which transforms surface integrals into volume integrals. This technique was already used as a main tool to homogenize the non homogeneous Neumann problem for the elliptic case by Cioranescu and Donato \cite{Ciora2}. 

Following  the technical lemmas contained in \cite[Section 3]{Ciora2}, let us introduce, for any $h\in L^{s'}(\partial F)$, $1\leq s'\leq \infty$, the linear form $\mu_{h}^\varepsilon$ on $W_0^{1,s}(\Omega)$ defined by
$$\langle \mu_{h}^\varepsilon,\varphi \rangle=\varepsilon \int_{\partial F_\varepsilon} h(x'/ \varepsilon,y_3)\varphi \,d\sigma(x')dy_3,\quad \forall \varphi\in W_0^{1,s}(\Omega),$$
with $1/s+1/s'=1$, such that
\begin{equation}\label{convergence_mu}
\mu_{h}^\varepsilon \to \mu_{h}\quad \text{strongly in }(W_0^{1,s}(\Omega))',
\end{equation}
where $\langle \mu_{h},\varphi \rangle=\mu_{h}\displaystyle\int_{\Omega}\varphi dx'dy_3$, with
$$\mu_{h}={1\over |Y'|}\int_{\partial F}h(y)d\sigma(y')dy_3.$$
In the particular case in which $h\in L^{\infty}(\partial F)$ or even when $h$ is constant, we have
\begin{equation*}
\mu_{h}^\varepsilon \to \mu_{h}\quad \text{strongly in }W^{-1,\infty}(\Omega).
\end{equation*}
In what follows, we shall denote by $\mu_1^\varepsilon$ the above introduced measure in the particular case in which $h=1$. Notice that in this case $\mu_{h}$ becomes $\mu_1=|\partial F'|/|Y'|$.

Observe that using \cite[Corollary 4.2]{Ciora3} with (\ref{continuity1}), we can deduce, for $\varepsilon \to 0$, 
\begin{equation*}
\varepsilon \int_{\partial F_\varepsilon} \gamma_0(\tilde u_\varepsilon(t))v d\sigma(x')dy_3=\langle \mu_{1}^\varepsilon,\widetilde \varPi_\varepsilon\tilde u_{\varepsilon {|_{\widetilde\Omega_\varepsilon}}}(t)v \rangle\to \mu_1\int_{\omega}u(t)v dx'={|\partial F'| \over |Y'|}\int_{\omega}u(t)v dx',
\end{equation*}
which integrating in time and using Lebesgue's Dominated Convergence Theorem, gives
$$\varepsilon \int_0^T{d \over dt}\varphi(t)\left(\int_{\partial F_\varepsilon} \gamma_0(\tilde u_\varepsilon(t))v d\sigma(x')dy_3\right) dt \to {|\partial F'| \over |Y'|}\int_0^T{d \over dt}\varphi(t)\left(\int_{\omega}u(t)v dx'\right)dt.$$
Moreover, using \cite[Corollary 4.2]{Ciora3} with (\ref{continuity1}), we can deduce, for $\varepsilon \to 0$, 
\begin{equation*}
\varepsilon\int_{\partial F_\varepsilon}
\gamma_{0}(u_\varepsilon(0))vd\sigma(x')dy_3=\langle \mu_{1}^\varepsilon,\widetilde \varPi_\varepsilon\tilde u_{\varepsilon {|_{\widetilde\Omega_\varepsilon}}}(0)v \rangle\to \mu_1\int_{\omega}u(0)v dx'={|\partial F'| \over |Y'|}\int_{\omega}u(0)v dx'.
\end{equation*}

On the other hand, note that using (\ref{convergence_mu}) with $s=2$, taking into account (\ref{hyp 0'}) and by Lebesgue's Dominated Convergence Theorem, we can deduce, for $\varepsilon \to 0$, 
$$\varepsilon \!\!\int_0^T \!\!\!\!\varphi(t)\!\!\int_{\partial F_\varepsilon}\!\!\!\! g(x',t)vd\sigma(x')dy_3dt=\int_0^T\!\!\!\!\varphi(t)\langle \mu_{1}^\varepsilon,g(x',t)v \rangle dt\to {|\partial F'| \over |Y'|}\int_0^T \!\!\!\!\varphi(t)\int_{\omega}g(x',t)vdx'dt.$$

{\bf Step 3}. 
Passing to the limit, as $\varepsilon\to 0$, in (\ref{system1}): all the terms in (\ref{system1}) pass to the limit, as $\varepsilon \to 0$, and therefore taking into account the previous steps, we get
 \begin{eqnarray*}
-\varphi(0)\left({|Y'_f|\over |Y'|}+{|\partial F'| \over |Y'|} \right)\int_{\omega} u(0)vdx'
-\left({|Y'_f|\over |Y'|}+{|\partial F'| \over |Y'|} \right)\int_0^T\dfrac{d}{dt}\varphi(t)\int_{\omega} u(t)vdx'dt\\[2ex]\nonumber
+\int_0^T\varphi(t)\int_{\Omega}\xi\cdot\nabla_{x'} vdx'dy_3dt+\kappa {|Y'_f|\over |Y'|}\int_0^T\varphi(t)\int_{\omega} u(t)v dx'dt\\[2ex]\nonumber
 ={|Y'_f|\over |Y'|}\int_0^T \varphi(t)\int_{\omega}  f(t)vdx'dt
 +{|\partial F'| \over |Y'|} \int_0^T \varphi(t)\int_{\omega}g(t)vdx'dt.\nonumber
 \end{eqnarray*}
 Hence, $\int_0^1\xi dy_3$ verifies
 \begin{equation}\label{equation_xi}
 \left({|Y'_f|\over |Y'|}\!+\!{|\partial F'| \over |Y'|} \right)\displaystyle{\partial u\over \partial t}\!-\!{\rm div}_{x'}\left(\int_0^1\xi dy_3\right)\!+\! {|Y'_f|\over |Y'|}\kappa u\!=\! {|Y'_f|\over |Y'|}f\!+\!{|\partial F'| \over |Y'|}g, 
 \end{equation}
 in $\omega\times (0,T).$
  
{\bf Step 4.} It remains now to identify $\int_0^1\xi dy_3$. We shall make use of the solutions of the cell problems (\ref{system_eta}). For any fixed $i=1,2$, let us define
\begin{equation}\label{definition_Psi}
\Psi_{i\varepsilon}(x')=\varepsilon\left(w_i\left({x'\over \varepsilon}\right)+y'_i \right)\quad \forall x'\in \omega_\varepsilon,
\end{equation}
where $y'=x'/\varepsilon$.

Recalling that $w_i$ is $Y'$-periodic, we obtain, in view of \cite[Theorem 2.6]{Ciora-Donato-book}, that
\begin{equation*}\label{convergence_Psi}
P_\varepsilon \Psi_{i\varepsilon}\rightharpoonup
x'_i \quad \textrm{weakly in} \quad H^1(\omega),
\end{equation*}
where $P_\varepsilon \Psi_{i\varepsilon}$ denotes the extension to $\omega$ given in \cite[Lemma 1]{Cioranescu}. Then, by Rellich-Kondrachov Theorem, we can deduce
\begin{equation}\label{convergence_Psi_fuerte}
P_\varepsilon \Psi_{i\varepsilon}\rightarrow
x'_i \quad \textrm{strongly in} \quad L^2(\omega).
\end{equation}
Let $\nabla_{x'} \Psi_{i\varepsilon}$ be the gradient of $\Psi_{i\varepsilon}$ in $\omega_\varepsilon$. Denote by $\widetilde{\nabla_{x'} \Psi_{i\varepsilon}}$ the extension by zero of $\nabla_{x'} \Psi_{i\varepsilon}$ inside the holes. From (\ref{definition_Psi}), we have
$$\widetilde{\nabla_{x'} \Psi_{i\varepsilon}}=\widetilde{\nabla_{y'}(w_i+y'_i)}=\widetilde{\nabla_{y'} w_i}(y')+e_i\chi_{Y'_f},$$
and taking into account \cite[Theorem 2.6]{Ciora-Donato-book}, we have
\begin{eqnarray}\label{limit_gradient_Psi}
\widetilde{\nabla_{x'} \Psi_{i\varepsilon}}\rightharpoonup
{1\over |Y'|} \int_{Y'_f}\left(e_i+\nabla_{y'} w_i(y') \right)dy' \quad \textrm{weakly in} \quad L^2(\omega).
\end{eqnarray}
On the other hand, it is not difficult to see that $ \Psi_{i\varepsilon}$ satisfies
\begin{equation}\label{system_Psi}
\left\{
\begin{array}{l}
\displaystyle -{\rm div}_{x'}\left(\nabla_{x'} \Psi_{i\varepsilon}\right)=0,   \text{\ in }\omega_\varepsilon,\\[2ex]
\nabla_{x'} \Psi_{i\varepsilon}\cdot \nu'_\varepsilon =0,  \text{\ on }\partial F'_\varepsilon.\end{array}
\right.
\end{equation}
Let $v\in \mathcal{D}(\Omega)$ with $v$ independent of $y_3$. Multiplying the first equation in (\ref{system_Psi}) by $v\tilde u_\varepsilon$ and integrating by parts over $\omega_\varepsilon$, we get
\begin{eqnarray}\label{igualdad_problema_Psi}
\int_{\omega_\varepsilon}\nabla_{x'} \Psi_{i\varepsilon}\cdot \nabla_{x'} v\,\tilde u_\varepsilon dx'+\int_{\omega_\varepsilon}\nabla_{x'} \Psi_{i\varepsilon}\cdot \nabla_{x'} \tilde u_\varepsilon vdx'=0.
\end{eqnarray}
On the other hand, we multiply system (\ref{PDE_dilatado}) by the test function $v \Psi_{i\varepsilon}$ and integrating by parts over $\widetilde \Omega_\varepsilon$, we obtain
\begin{eqnarray}\label{formula_indetificaXI}
\dfrac{d}{dt}\left(\int_{\Omega}\chi_{\widetilde\Omega_\varepsilon}\widetilde \varPi_\varepsilon\tilde u_\varepsilon v P_\varepsilon\Psi_{i\varepsilon}dx'dy_3\right)+\varepsilon\,\dfrac{d}{dt}\left(\int_{\partial F_\varepsilon}
\gamma_{0}(\tilde u_\varepsilon)v\gamma_0(\Psi_{i\varepsilon})d\sigma(x')dy_3\right)\nonumber\\[2ex]
+\int_{\widetilde\Omega_\varepsilon} \nabla_{x'} \tilde u_\varepsilon \cdot \nabla_{x'} v \Psi_{i\varepsilon}dx'dy_3+\int_{\widetilde\Omega_\varepsilon} \nabla_{x'} \tilde u_\varepsilon \cdot \nabla_{x'} \Psi_{i\varepsilon}vdx'dy_3
+\kappa \int_{\Omega}\chi_{\widetilde\Omega_\varepsilon}\widetilde \varPi_\varepsilon\tilde u_\varepsilon v P_\varepsilon \Psi_{i\varepsilon} dx'dy_3
\\[2ex]
=\int_{\Omega}\chi_{\widetilde\Omega_\varepsilon}f(x',t) v P_\varepsilon \Psi_{i\varepsilon} dx'dy_3
 +\varepsilon\, \int_{\partial F_\varepsilon}
g(x',t)v\gamma_0(\Psi_{i\varepsilon})d\sigma(x')dy_3, \nonumber
 \end{eqnarray}
 in $\mathcal{D}'(0,T)$. 
 
Using (\ref{igualdad_problema_Psi}) in (\ref{formula_indetificaXI}), we have
\begin{eqnarray*}
\dfrac{d}{dt}\left(\int_{\Omega}\chi_{\widetilde\Omega_\varepsilon}\widetilde \varPi_\varepsilon\tilde u_\varepsilon v P_\varepsilon\Psi_{i\varepsilon}dx'dy_3\right)+\varepsilon\,\dfrac{d}{dt}\left(\int_{\partial F_\varepsilon}
\gamma_{0}(\tilde u_\varepsilon)v\gamma_0(\Psi_{i\varepsilon})d\sigma(x')dy_3\right)\nonumber\\[2ex]
+\int_{\Omega} \tilde \xi_\varepsilon \cdot \nabla_{x'} v P_\varepsilon\Psi_{i\varepsilon}dx'dy_3-\int_{\Omega}\widetilde{\nabla_{x'} \Psi_{i\varepsilon}}\cdot \nabla_{x'} v\,\widetilde \varPi_\varepsilon\tilde u_\varepsilon dx'dy_3
+\kappa \int_{\Omega}\chi_{\widetilde\Omega_\varepsilon}\widetilde \varPi_\varepsilon\tilde u_\varepsilon v P_\varepsilon \Psi_{i\varepsilon} dx'dy_3
\\[2ex]
=\int_{\Omega}\chi_{\widetilde\Omega_\varepsilon}f(x',t) v P_\varepsilon \Psi_{i\varepsilon} dx'dy_3
 +\varepsilon\, \int_{\partial F_\varepsilon}
g(x',t)v\gamma_0(\Psi_{i\varepsilon})d\sigma(x')dy_3, \nonumber
 \end{eqnarray*}
 in $\mathcal{D}'(0,T)$. 
 
We consider $\varphi\in C_c^1([0,T])$ such that $\varphi(T)=0$ and $\varphi(0)\ne 0$. Multiplying by $\varphi$ and integrating between $0$ and $T$, we have
\begin{eqnarray}\label{New_formula}
-\varphi(0)\left(\int_{\Omega}\chi_{\widetilde\Omega_\varepsilon}\widetilde \varPi_\varepsilon\tilde u_\varepsilon(0) v P_\varepsilon\Psi_{i\varepsilon}dx'dy_3\right)-\int_0^T \dfrac{d}{dt}\varphi(t)\left(\int_{\Omega}\chi_{\widetilde\Omega_\varepsilon}\widetilde \varPi_\varepsilon\tilde u_\varepsilon(t) v P_\varepsilon\Psi_{i\varepsilon}dx'dy_3\right)dt\nonumber\\[2ex]
-\varepsilon\,\varphi(0)\left(\int_{\partial F_\varepsilon}
\gamma_{0}(\tilde u_\varepsilon(0))v\gamma_0(\Psi_{i\varepsilon})d\sigma(x')dy_3\right)-\varepsilon\,\int_0^T\dfrac{d}{dt}\varphi(t)\left(\int_{\partial F_\varepsilon}
\gamma_{0}(\tilde u_\varepsilon(t))v\gamma_0(\Psi_{i\varepsilon})d\sigma(x')dy_3\right)dt\nonumber\\[2ex]
+\int_0^T\varphi(t)\int_{\Omega} \tilde \xi_\varepsilon \cdot \nabla_{x'} v P_\varepsilon\Psi_{i\varepsilon}dx'dy_3dt-\int_0^T \varphi(t)\int_{\Omega}\widetilde{\nabla_{x'} \Psi_{i\varepsilon}}\cdot \nabla_{x'} v\,\widetilde \varPi_\varepsilon\tilde u_\varepsilon dx'dy_3dt\\[2ex]
+\kappa \int_0^T \varphi(t)\int_{\Omega}\chi_{\widetilde\Omega_\varepsilon}\widetilde \varPi_\varepsilon\tilde u_\varepsilon(t) v P_\varepsilon \Psi_{i\varepsilon} dx'dy_3dt\nonumber
\\[2ex]
=\int_0^T \varphi(t)\int_{\Omega}\chi_{\widetilde\Omega_\varepsilon}f(x',t) v P_\varepsilon \Psi_{i\varepsilon} dx'dy_3dt
 +\varepsilon\, \int_0^T \varphi(t)\int_{\partial F_\varepsilon}
g(x',t)v\gamma_0(\Psi_{i\varepsilon})d\sigma(x')dy_3dt. \nonumber
 \end{eqnarray}

Now, we have to pass to the limit, as $\varepsilon\to 0$. We will focus on the terms which involve the gradient. Taking into account (\ref{convergence_Psi_fuerte}), we reason as in steps 1 and 2 for the others terms. 
 
 Firstly, using (\ref{converge_gradiente}) and (\ref{convergence_Psi_fuerte}), we have 
\begin{eqnarray*}\label{limite_1}
\int_0^T \varphi(t)\int_{\Omega} \tilde \xi_\varepsilon \cdot \nabla_{x'} v P_\varepsilon \Psi_{i\varepsilon}dx'dy_3dt\to \int_0^T \varphi(t)\int_{\Omega}\xi \cdot \nabla_{x'} v\, x'_idx'dy_3dt,
\end{eqnarray*}
and by (\ref{converge_initial_data}), (\ref{limit_gradient_Psi}) and Lebesgue's Dominated Convergence Theorem, we obtain
\begin{eqnarray*}\label{limite_2}
\int_0^T\varphi(t)\int_{\Omega}\widetilde{\nabla_{x'} \Psi_{i\varepsilon}}\cdot \nabla_{x'} v\,\widetilde \varPi_\varepsilon\tilde u_\varepsilon dx'dy_3dt\to {1\over |Y'|}\int_0^T \varphi(t)\int_{\omega}\left( \int_{Y'_f}\left(e_i+\nabla_{y'} w_i \right)dy'\right)\cdot \nabla_{x'} v\, u\,dx'dt.
\end{eqnarray*}
Therefore, when we pass to the limit in (\ref{New_formula}), we obtain
\begin{eqnarray*}
-\varphi(0)\left({|Y'_f|\over |Y'|}+{|\partial F'| \over |Y'|} \right)\left(\int_{\omega} u(0)vx'_idx'\right)
-\left({|Y'_f|\over |Y'|}+{|\partial F'| \over |Y'|} \right)\int_0^T\dfrac{d}{dt}\varphi(t)\int_{\omega} u(t)vx'_idx'dt\\[2ex]
+\int_0^T \varphi(t)\int_{\Omega}\xi \cdot \nabla_{x'} v\, x'_idx'dy_3dt-{1\over |Y'|}\int_0^T \varphi(t)\int_{\omega}\left( \int_{Y'_f}\left(e_i+\nabla_{y'} w_i \right)dy'\right)\cdot \nabla_{x'} v\, u(t)\,dx'dt\\[2ex]\nonumber
+\kappa {|Y'_f|\over |Y'|}\int_0^T\varphi(t)\int_{\omega} u(t)vx'_i dx'dt\\[2ex]
={|Y'_f|\over |Y'|}\int_0^T \varphi(t)\int_{\omega}  f(t)vx'_idx'dt
 +{|\partial F'| \over |Y'|} \int_0^T \varphi(t)\int_{\omega}g(t)vx'_idx'dt.\\[2ex]
 \end{eqnarray*}
 Using Green's formula and equation (\ref{equation_xi}), we have
\begin{eqnarray*}
-\int_0^T \varphi(t)\int_{\Omega}\xi \cdot \nabla_{x'} x'_i\, vdx'dy_3dt+{1\over |Y'|}\int_0^T \varphi(t)\int_{\omega}\left( \int_{Y'_f}\left(e_i+\nabla_{y'} w_i \right)dy'\right)\cdot \nabla_{x'} u\,vdx'dt=0.
\end{eqnarray*}
The above equality holds true for any $v\in \mathcal{D}(\Omega)$ independent of $y_3$, and $\varphi\in C_c^1([0,T])$. This implies that
 \begin{eqnarray*}
 -\left(\int_0^1\xi dy_3\right) \cdot \nabla_{x'} x'_i +{1\over |Y'|}\left( \int_{Y'_f}\left(e_i+\nabla_{y'} w_i \right)dy'\right)\cdot \nabla_{x'} u
 =0,\quad \text{ in }\omega\times (0,T).
\end{eqnarray*}
We conclude that
 \begin{eqnarray}\label{identificacion_xi}
{\rm div}_{x'}\left(\int_0^1\xi dy_3\right)= {\rm div}_{x'}\left(Q \nabla_{x'} u\right),
\end{eqnarray}
where $Q=((q_{ij}))$, $1\leq i,j\leq 2$, is given by 
\begin{eqnarray*}\label{A_limite}
q_{ij}={1\over |Y'|} \int_{Y'_f}\left(e_{i}+\nabla_{y'} w_i \right)\cdot e_j \,dy'.
\end{eqnarray*}
Observe that if we multiply system (\ref{system_eta}) by the test function $w_j$, integrating by parts over $Y'_f$, we obtain
\begin{eqnarray*}
\int_{Y'_f}(e_i+\nabla_{y'} w_i)\cdot \nabla_{y'} w_j dy'=0,
\end{eqnarray*}
then we conclude that $q_{ij}$ is given by (\ref{matrix}).

{\bf Step 5.} Finally, thanks to (\ref{equation_xi}) and (\ref{identificacion_xi}), we observe that $u$ satisfies the first equation in (\ref{limit_problem}). A weak solution of (\ref{limit_problem}) is any function $u$, satisfying
\begin{equation*}
u\in {C}([0,T];L^{2}\left(  \omega\right)  ),\quad \text{for all }T>0,
\end{equation*}
\begin{equation*}
u\in L^{2}(0,T;H_{0}^{1}\left(  \omega\right)  ),\quad \text{for all }T>0,
\end{equation*}
\begin{equation*}
\displaystyle \left({|Y'_f|\over |Y'|}+{|\partial F'| \over |Y'|} \right) \dfrac{d}{dt}(u(t),v)+(Q\nabla_{x'} u(t),\nabla_{x'} v)+{|Y'_f|\over |Y'|}\kappa(u(t),v)= {|Y'_f|\over |Y'|}(f(t),v)+{|\partial F'| \over |Y'|}(g(t),v),\quad \text{in  }\mathcal{D}'(0,T),
\end{equation*}
for all $v\in H_0^1(\omega)$, and
\begin{equation*}
 u(0)=u_0.
\end{equation*}

Since the homogenized matrix $Q$ is positive-definite (see \cite[Theorem 4.7]{Ciora2}), applying a slight modification of \cite[Chapter 2,Theorem 1.4]{Lions}, we obtain that the problem (\ref{limit_problem}) has a unique solution, and therefore Theorem \ref{Main} is proved.

%%%%%%%%%%%%%%%%%%%%%%%%%%%FIN

\end{document}